\documentclass{amsart}

\usepackage[left=3cm,top=3cm,right=3cm,bottom=3cm]{geometry}
\usepackage{amssymb, amsmath, amsthm}
\usepackage{verbatim}
\usepackage{comment}
\usepackage{graphicx}
\usepackage{enumerate}
\usepackage{color}

\DeclareGraphicsExtensions{.pdf,.jpg,.png}

\theoremstyle{plain}
\newtheorem{theorem}{Theorem}

\newtheorem{conjecture}[theorem]{Conjecture}
\newtheorem{problem}[theorem]{Problem}

\theoremstyle{definition}

\newcommand{\A}{\mathcal{A}}

\makeatletter

\begin{document}

\title{On sums and products along the edges, II}

\author{Noga Alon}
\address{\noindent  Princeton University, Princeton, NJ 08544 and
Tel Aviv University, Tel Aviv 69978,
Israel} 
\email{nalon@math.princeton.edu}
\thanks{Research supported in part by NSF grant DMS-2154082.}
\author{Imre Ruzsa}
\address{\noindent Alfr\'ed R\'enyi Institute of Mathematics, 
H-1364 Budapest, Hungary}
\email{ruzsa@renyi.hu}
\thanks{Research supported in part by an OTKA NK 133819 grant}
\author{J\'{o}zsef Solymosi}
\address{\noindent Department of Mathematics, 
University of British Columbia, 1984 Mathematics Road, 
Vancouver, BC, V6T 1Z2, Canada and Obuda University, H-1034 Budapest, Hungary }
\email{solymosi@math.ubc.ca}
\thanks{Research supported in part by an NSERC and an OTKA NK 133819 grant}

\maketitle

\begin{abstract}
This note is a continuation of an earlier paper of the authors
\cite{ARS}. We describe improved constructions addressing a question
of Erd\H{o}s and Szemer\'edi on sums and products of real numbers 
along the edges of a graph. 
We also add a few observations about related versions
of the problem.
\end{abstract}

\section{Introduction}

In this note, we describe an improved construction addressing
a question of Erd\H{o}s
and Szemer\'edi about sums and products along the edges of a graph. We also
mention some related problems. The main
improvement is obtained by a simple modification of the
construction  in \cite{ARS} which works for real numbers,
instead of the integers considered there.

In their original paper  Erd\H{o}s and Szemer\'edi \cite{ERDSZE}
considered sum and product along the edges of graphs.
Let $G_n$ be a graph on $n$ vertices, $v_1, v_2, \ldots, v_n,$ with
$n^{1+c}$ edges for some real $c>0$. Let $\A$ be an $n$-element
set of real numbers, $\A=\{a_1, a_2,\ldots, a_n\}.$ The {\em sumset of
$\A$ along $G_n$}, denoted by $\A+_{G_n}\A$, is the set $\{a_i+a_j |
(i,j)\in E(G_n)\}.$ The product set along $G_n$ is defined similarly,
$$\A\cdot_{G_n}\A=\{a_i\cdot a_j | (i,j)\in E(G_n)\}.$$

The Strong Erd\H{o}s-Szemer\'edi Conjecture, which 
was refuted in \cite{ARS}, is the following. 

\begin{conjecture}\label{Strong}\cite{ERDSZE}
For every $c>$ and $\varepsilon >0,$ there is a threshold, $n_0,$
such that if $n\geq n_0$ then for any $n$-element subset of reals
$\A\subset \mathbb{R}$ 
and any graph $G_n$ with $n$ vertices and at least $n^{1+c}$ edges
$$
|\A+_{G_n}\A|+|\A\cdot_{G_n}\A|\geq |\A|^{1+c-\varepsilon}.
$$
\end{conjecture}

Now the question is to find dense graphs with small sumset and
product set along the edges. Here we extend the construction
in \cite{ARS}. The improvement follows by considering real numbers,
instead of integers only.

\section{Constructions}
\subsection{Sum-product along edges with real numbers}
Here we extend our earlier construction so that we get 
better bounds in a range of edge densities.
In our previous paper for arbitrary large $m_0,$ 
we constructed a set of integers, $\A,$ and a  graph
on $|\A|=m\geq m_0$ vertices, $G_m,$ with $\Omega(m^{5/3}/\log^{1/3}m)$
edges such that
$$
|\A+_{G_m}\A|+|\A\cdot_{G_m}\A|= O\left((|\A|\log |\A|)^{4/3}\right).
$$

Thus we had a graph on $m$ vertices and roughly $m^{2-c}$ 
edges with roughly $m^{2-2c}$ sums and products along the edges 
for $c=1/3.$ 
In the following construction, we show a similar bound in a range 
covering all $0 \leq c \leq 2/5$.
In what follows, it is convenient to ignore the
logarithmic terms. We thus use now the common notation
$f=\tilde{O}(g)$ for two functions $f(n)$ and $g(n)$ 
to denote that there are absolute positive constants
$c_1,c_2$ so that 
$ 
f(n) \leq c_1 g(n) (\log g(n))^{c_2}
$
for all admissible values of $n$.  The notation $f = \tilde{\Omega}(g)$
means that $g=\tilde{O}(f)$ and $f=\tilde{\Theta}(g)$ denotes that
$f = \tilde{\Omega}(g)$ and $g=\tilde{O}(f)$.
\medskip

\begin{theorem}
\label{sumprod}
For arbitrary large $m_0,$ 
and parameter $\alpha,$ where $0\leq\alpha\leq 1/5,$ 
there is a set of reals, $\A,$ and a graph
on $|\A|=m\geq m_0$ vertices, $G_m,$ with 
$$
\tilde{\Omega}\left( m^{2-2\alpha}\right)
$$
edges such that
$$
|\A+_{G_m}\A|+|\A\cdot_{G_m}\A|= 
\tilde{O}\left(|\A|^{2-4\alpha}\right).
$$
\end{theorem}

{\em Proof:} \,
It is easier to describe the construction using prime numbers only. 
We get a slightly larger exponent in the hidden logarithmic factor, 
but we are anyway ignoring these factors here. 
The set of primes is denoted by $\mathbb{P}$ here.
We define the set $\A$ first and then the graph using the 
parameter  $\alpha.$

$$ 
\A:=\left\{\frac{u\sqrt{w}}{\sqrt{v}} \text { } | 
\text { } u,v,w\in\mathbb{P}
\text{ distinct and }  v,w
\leq n^{\alpha}, u\leq n^{1-2\alpha} \right\}. 
$$

It is clear that distinct choices
of $3$-tuples $u,v,w$ lead to distinct reals. Thus with this
choice of parameters, the size of $\A$ is
$\tilde{\Theta}(n).$
We are going to define a graph $G_{m}$ with
vertex set $\A,$ where $|\A|=m= \tilde{\Theta}\left(n\right).$ 
Two elements, $a,b\in\A$
are connected by an edge if in the definition of $\A$ above 
$a=\frac{u\sqrt{w}}{\sqrt{v}}$
and $b=\frac{z\sqrt{v}}{\sqrt{w}}.$ 
Since the degree of every vertex here is
$\tilde{\Theta}(n^{1-2\alpha})$
the number of edges is
$$
\tilde{\Omega}\left(m^{2-2\alpha}\right).
$$
The products of pairs of elements of $\A$ along an edge of
$G_{m}$ are integers of size at most 
$$n^{2-4\alpha}=
\tilde{O}\left(m^{2-4\alpha}\right) .
$$ 
The sums
along the edges are  of the form

$$\frac{u\sqrt{w}}{\sqrt{v}}+\frac{z\sqrt{v}}{\sqrt{w}}
=\frac{wu+vz}{\sqrt{vw}}.$$

The number of possibilities for the denominator
is at most $n^{2\alpha}$ and
the numerator is a positive integer of size at 
most $2n^{1-\alpha}$, hence
the number of sums is, at most
$$
O(n^{1+\alpha})
=\tilde{O}\left(m^{2-(1-\alpha)}\right).
$$  

The sum is asymptotically smaller than the product set, as long as $1-\alpha > 4\alpha,$ i.e. $\alpha < 1/5.$

\qed

\medskip

Based on this construction, one can easily get examples of sparser 
graphs, simply taking smaller copies of $G_m$ and leaving 
other vertices isolated. 

\begin{theorem}
For every parameters $0 \leq \nu\leq 3/5$  and 
$n_0$ there are $n>n_0$,
an $n$-element set of reals, 
$\A\subset
\mathbb{R},$ and a graph $H_n$ with $\tilde{\Omega}
(n^{1+\nu})$ edges such that
$$|\A+_{H_n}\A|+|\A\cdot_{H_n}\A|=\tilde{O}
\left(|\A|^{3(1+\nu)/4 }\right).$$
\end{theorem}

{\em Proof:} \, The construction of Theorem \ref{sumprod} with
$\alpha=1/5$ supplies a set of $m$ reals and a graph with
$\tilde{\Omega}(m^{8/5})$ edges so that the number of
sums and products along the edges is at most $\tilde{O}(m^{6/5})$.
Take this construction with $m=n^{5(1+\nu)/8} (\leq n)$ and add
to it $n-m$ isolated vertices assigning to them arbitrary distinct
reals that differ from the ones used already.
\qed

A similar statement holds for integers too.

\begin{theorem}
\label{t4}
For every parameters $0 \leq \nu\leq 2/3$ and $n_0$  there 
are $n>n_0$, an $n$-element set of integers 
$\A$, and a graph $H_n$ with $\tilde{\Omega}(n^{1+\nu})$ edges such that
$$|\A+_{H_n}\A|+|\A\cdot_{H_n}\A|=
\tilde{O}\left(|\A|^{4(1+\nu)/5 }\right).$$
\end{theorem}

This follows as in the real case by starting with the construction of 
\cite{ARS} that gives a set of $m$ integers and a graph with
$\tilde{\Omega}(m^{5/3})$ edges  so that the number of sums and
products along the edges is at most $\tilde{O}(m^{4/3})$. 
This construction with $m=n^{3(1+\nu)/5} \leq n$ together with
$n-m$ isolated vertices with arbitrary $n-m$ new integers implies
the statement above.

\subsection{Matchings}

A particular variant of the sum-product problem for integers 
is the following:
\begin{problem}
Given two $n$-element sets of integers, $A=\{a_1,\ldots ,a_n\}$
and $B=\{b_1,\ldots ,b_n\} $ let us define a sumset and a product set as 
$$S=\{a_i+b_i |  1\leq i\leq n\}\text{  }
and\text{  }
P=\{a_i\cdot b_i |  1\leq i\leq n\}.$$
Erd\H{o}s and Szemer\'edi conjectured that 
\begin{equation}\label{SZP}
|P|+|S|=\Omega(n^{1/2+c})
\end{equation}
for some constant $c>0.$ 
\end{problem} 

The best-known lower bound is due to Chang \cite{MCC}, 
who proved that \[|P|+|S|\geq n^{1/2}\log^{1/48} n.\]

It was shown recently in \cite{SS} that under the assumption of a special
case of the Bombieri-Lang conjecture \cite{CHM}, one can take $c=1/10$ in
equation (\ref{SZP}), i.e. $|P|+|S|=\Omega(n^{3/5}),$ even for multisets.
\medskip

\begin{theorem}\label{matching}\cite{SS}
Let $M=\{(a_i,b_i) | 1\leq i\leq
n\}$ be a set of distinct pairs of integers.
If $P$ and $S$ are defined as above, then under the hypothesis of
the Bombieri-Lang conjecture $|P|+|S|=\Omega(n^{1/2+c})$  with $c=1/10.$
\end{theorem}

\medskip
If multisets are allowed, and the only requirement is that the pairs
assigned to distinct edges of the matching are distinct, then any
construction of a graph with $n$ edges yields a construction of a
matching of size $n$. It thus follows from 
\cite[Theorem 3 ]{ARS} (or from Theorem \ref{t4} here) 
that for the multiset version
there is, for arbitrarily large $n$, an example of a matching $M$ of
size $n$ as above, with $n$ distinct pairs of integers $(a_i,b_i)$,
so that $|P|+|S|=\tilde{O}(n^{4/5}).$ This shows
that the statement of Theorem \ref{matching} cannot be 
improved beyond an extra $1/5$ in the exponent.
\medskip

\section{Lower bounds}
In \cite{ARS}, we followed Elekes' method using point-line incidence bounds to give a lower bound on the sum-product problem along the edges of a graph. For sparser graphs, Oliver Roche-Newton improved our bound, extending the range where a non-trivial bound can be established. He proved the following

\begin{theorem}[Theorem 6.1 in \cite{ORN}]\label{Oliver}
    For arbitrary set of reals, $\A,$ and a graph
on $|\A|=m$ vertices, $G_m,$ with 
$$
\tilde{\Omega}\left( m^{2-2\alpha}\right)
$$
edges the following bound holds:
$$
|\A+_{G_m}\A|+|\A\cdot_{G_m}\A|= 
\tilde{\Omega}\left(|\A|^{\frac{9-12\alpha}{8}}\right).
$$
\end{theorem}

The result follows from applying an Elekes-Szab\'o type bound on the intersection size of polynomials and Cartesian products. Roche-Newton used the bound from \cite{RSZ}, however, a better result follows from the recent improvement in \cite{SZ}. 

\begin{theorem}\label{ESZ}[Theorem 1.4 in \cite{SZ}]\label{SoZh}
Let $f\in {\mathbb{C}}[x,y,z]$ be an irreducible polynomial. Then at least one of the following is true.
\begin{itemize}

	\item[(A)] For all finite sets $A,B,C\subset{\mathbb{R}}$ with $|A|\leq |B|\leq |C|$, we have
	\begin{equation*}
	| (A\times B \times C)\cap Z(f)|=\tilde{O} (|A| |B| |C|)^{4/7} + |B| |C|^{1/2},
	\end{equation*}
	where the implicit constant depends on the degree of $f$.

	\item[(B)] After possibly permuting the coordinates $x,y,z$, we have $f(x,y,z) = g(x,y)$, for some bivariate polynomial $g$. 

	\item[(C)] $f$ encodes additive group structure.{\footnote{When $f(x,y,z)$ is of the special form $h(x,y)-z$, then $f$ encodes additive structure if and only if $h$ has the form $h(x,y) = p(q(x)+r(y))$ or $h(x,y) = p(q(x)r(y))$ for univariate polynomials $p,q,r$.}}
\end{itemize}
\end{theorem}

Now we state a new lower bound on the size of the sumset and product set along the edges of a graph.

\begin{theorem}
For arbitrary set of reals, $\A,$ and a graph
on $|\A|=m$ vertices, $G_m,$ with 
$$
\tilde{\Omega}\left( m^{2-2\alpha}\right)
$$
edges the following bound holds:
$$
|\A+_{G_m}\A|+|\A\cdot_{G_m}\A|= 
\tilde{\Omega}\left(|\A|^{\frac{5-7\alpha}{4}}\right).
$$
\end{theorem}

{\em Proof:} \,
For the proof we can follow the arguments in \cite{ORN} and use the new Elekes-Szab\'o type bound from Theorem \ref{ESZ}. We consider the zero set of the polynomial
\[
f(x,y,z)=x(y-x)-z,
\]
and its intersection with the Cartesian product $\A\times\{\A+_{G_m}\A\}\times \{\A\cdot_{G_m}\A\}$. 
Every edge in $G_m$ which connects vertices $a$ and $b$ determines an intersection point, by $x=a,$ $y=a+b$ and $z=ab$.
This is the polynomial variant of Elekes' original sum-product bound in \cite{E} where he considered lines $\alpha(X-\beta)-Y=0$ with $\alpha,\beta\in \A$ and $X\in \A+\A, Y\in \A\A.$ 
As it was shown in \cite{ORN}, for this polynomial Part A applies from Theorem \ref{SoZh}. From that, we have the bound 

\[
m^{2-2\alpha} = \tilde{O} \left( \left(|\A| |\A+_{G_m}\A| |\A\cdot_{G_m}\A|\right)^{4/7} + |\A+_{G_m}\A| |\A\cdot_{G_m}\A|^{1/2}\right)
\]

which implies 
\[
|\A+_{G_m}\A|+|\A\cdot_{G_m}\A|= 
\tilde{\Omega}\left(|\A|^{\frac{5-7\alpha}{4}}\right).
\]

\qed

\section{Remarks}
There is still a gap between the lower bound and our construction. It is inevitable as long as the original sum-product conjecture is open. Our construction goes to the conjectured optimum as the graph is getting denser. The lower bound approaches Elekes' bound \cite{E}.

\begin{figure}[h]

\centering
\includegraphics[width=0.4\textwidth]{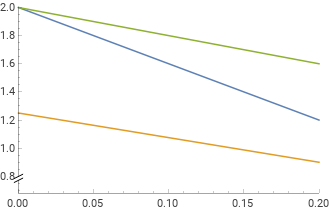}
\caption{The exponents in the upper and lower bounds when the number of edges is $m^{2-2\alpha}$ (top line) and $0<\alpha<1/5$}
\end{figure}


\begin{thebibliography}{00}

\bibitem{ARS} 
N. Alon, I. Ruzsa and J. Solymosi,  
{\em Sums, products, and ratios along the edges of a graph}, 
Publ. Mat. 64 (2020), no. 1, 143--155.

\bibitem{CHM} L. Caporaso, J. Harris, and B. Mazur, {\em Uniformity of
rational points}, J. Amer. Math.  Soc. 10 (1997), no. 1, 1--35.

\bibitem{MCC} 
{\sc M.-C. Chang, }
{\em  On problems of Erd\H{o}s and Rudin, }
Journal of Functional Analysis
Volume 207, Issue 2, (2004), 444--460.

\bibitem{E}
Gy. Elekes, {\em On the number of sums and products}
Acta Arithmetica, LXXXI.4, (1997) 365--367.


\bibitem{ERDSZE} P. Erd\H{o}s and E. Szemer\'{e}di, {\em On sums and
products of integers}, Studies in pure mathematics, Birkh\"{a}user,
Basel (1983) 213--218.

\bibitem{ORN} 
O. Roche-Newton,
SIAM J. Discrete Math.  Vol. 35, No. 1, pp. 194--204.

\bibitem{RSZ} O.~E.~Raz, M. Sharir, and F. de Zeeuw. {\em Polynomials vanishing on Cartesian products: The Elekes-Szab\'o theorem revisited} {Duke Math. J.} 165(18):3517--3566, 2016. 

\bibitem{RUD} W. Rudin, 
{\em Trigonometric series with gaps,}
J. Math. Mech. 9 (1960) 203--227.

\bibitem{SS} I.D. Shkredov and J. Solymosi,
{\em The Uniformity Conjecture in Additive Combinatorics}
SIAM Journal on Discrete Mathematics  Vol. 35, Iss. 1 (2021) 307--321.

\bibitem{SZ} J. Solymosi and J. Zahl,
{\em Improved Elekes-Szab\'o type estimates using proximity}, 
Journal of Combinatorial Theory, Series A. (to appear)
arXiv:2211.13294 [math.CO]

\end{thebibliography}
\end{document}